# Mixing property and pseudo random sequences

**Makoto Mori**[1]

*Nihon University*

**Abstract:** We will give a summary about the relations between the spectra of the Perron–Frobenius operator and pseudo random sequences for 1-dimensional cases.

There are many difficulties to construct general theory of higher-dimensional cases. We will give several examples for these cases.

## 1. Perron–Frobenius operator on one-dimensional dynamical systems

Let $I = [0, 1]$ and $F$ be a piecewise monotonic and $C^2$ transformations on it. We assume $F$ is expanding in the following sense:

$$\xi = \liminf_{n\to\infty} \operatorname*{ess\,inf}_{x\in I} \frac{1}{n} \log |F^{n\prime}(x)| > 0.$$

**Definition 1.1.** An operator $P\colon L^1 \to L^1$ defined by

$$Pf(x) = \sum_{y\,:\,F(y)=x} f(y)|F'(y)|^{-1}$$

is called the Perron–Frobenius operator associated with $F$.

It is well known that the spectra of the Perron–Frobenius operator $P$ associated with the transformation $F$ determine the ergodic properties of the dynamical system $(I, F)$. Roughly speaking:

1. The dimension of the eigenspace associated with eigenvalue 1 equals the number of the ergodic components.
2. The eigenfunction $\rho$ of $P$ associated with the eigenvalue 1 for which $\rho \geq 0$ and $\int \rho \, dx = 1$ becomes the density function of the invariant probability measure.
3. If 1 is simple and there exists no eigenvalue modulus 1 except 1, then the dynamical system is mixing.
4. The second greatest eigenvalue determines the decay rate of correlation.

To be more precise, the Perron–Frobenius operator satisfies:

1. $P$ is contracting, and 1 is an eigenvalue. Hence, its spectrum radius equals 1.
2. $P$ is positive. Hence there exists a nonnegative eigenfucntion associated with the eigenvalue 1.

---

[1]Department of Mathematics, College of Humanities and Sciences, Nihon University, Japan, e-mail: mori@math.chs.nihon-u.ac.jp

3. On th unit circle, the eigenvalues of $P$ coincide with the eigenvalues of the unitary operator $Uf(x) = f(F(x))$.
4. The essential spectrum radius of $P\colon L^1 \to L^1$ equals 1, that is, it is of no use to consider decay rate for the functions in $L^1$.
5. When we restrict $P$ to $BV$ the set of functions with bounded variation, then its essential spectrum radius equals $e^{-\xi}$. Hence, the decay rate for a function with bounded variation is at most $e^{-\xi}$.

Hereafter, we only consider $P$ restricted to $BV$. Since $P$ is not compact, it is not easy to calculate its spectra. Hofbauer and Keller ([1]) proved that they coincide with the singularities of the zeta function

$$\zeta(z) = \exp\left[\sum_{n=1}^{\infty} \frac{z^n}{n} \sum_{p\colon F^n(p)=p} |F^{n'}(p)|^{-1}\right].$$

The zeta function has the radius of convergence 1, and it has a meromorphic extension to $e^\xi$. But it is not easy to study the singularities in $|z| > 1$.

Since $F$ is piecewise monotone, there exists a finite set $\mathcal{A}$ which we call alphabet with the following properties.

1. For each $a \in \mathcal{A}$, an interval $\langle a \rangle$ corresponds, and $\{\langle a \rangle\}_{a \in \mathcal{A}}$ forms a partition of $I$,
2. $F$ is monotone on each $\langle a \rangle$, and it has a $C^2$ extension to the closure of $\langle a \rangle$.

We define
$$\operatorname{sgn} a = \begin{cases} +1 & \text{if } F \text{ is monotone increasing in } \langle a \rangle, \\ -1 & \text{otherwise.} \end{cases}$$

We call $F$ Markov if there exists an alphabet $\mathcal{A}$, such that, if $F(\langle a \rangle)^o \cap \langle b \rangle \neq \emptyset$ then $F(\langle a \rangle) \supset \langle b \rangle^o$, where $J^o$ denotes the inner of a set $J$.

We will summarize notations which we use in this article.

1. A sequence of symbols $a_1 \cdots a_n$ ($a_i \in \mathcal{A}$) is called a word and
   (a) $|w| = n$ (the length of a word),
   (b) $\langle w \rangle = \bigcap_{i=1}^{n} F^{-i+1}(\langle a_i \rangle)$ (the interval corresponding to a word),
   (c) We denote the empty word by $\epsilon$, and define $|\epsilon| = 0$, $\operatorname{sgn} \epsilon = +1$ and $\langle \epsilon \rangle = I$,
   (d) we call $w$ is admissible if $\langle w \rangle \neq \emptyset$. We denote the set of admissible words with length $n$ by $\mathcal{W}_n$ and $\mathcal{W} = \cup_{n=0}^{\infty} \mathcal{W}_n$,
   (e) $\operatorname{sgn} w = \prod_{i=1}^{n} \operatorname{sgn} a_i$.
2. For a point $x \in I$, $a_1^x a_2^x \cdots$ is called the expansion of $x$ which is defined by
$$F^n(x) \in \langle a_{n+1} \rangle.$$
   We usually identify $x$ and its expansion.
3. For a sequence of symbols $s = a_1 a_2 \cdots$ ($a_i \in \mathcal{A}$),
   (a) $s[n,m] = a_n a_{n+1} \cdots a_m$ ($n \leq m$),
   (b) $s[n] = a_n$.

4. For a word $w = a_1 \cdots a_n$ and $x \in I$, $wx$ is an infinite sequence of symbols defined by
$$wx = a_1 \cdots a_n a_1^x a_2^x \cdots.$$

We call $wx$ exists if there exists a point $y$ whose expansion equals $wx$.

We have a natural order on $\mathcal{A}$, that is, $a < b$ when $x < y$ for any $x \in \langle a \rangle$ and $y \in \langle b \rangle$. Then we define an order on words by

1. $w < w'$ if $|w| < |w'|$,
2. $w = a_1 \cdots a_n$, $w' = b_1 \cdots b_n$ and if $a_i = b_i$ for $k < i \leq n$ and $a_k < b_k$, then
$$\begin{aligned} w < w' & \quad \text{if sgn } a_{k+1} \cdots a_n = +1, \\ w > w' & \quad \text{otherwise.} \end{aligned}$$

## 2. Pseudo random sequences

We call a sequence $\{x_n\}_{n=1}^\infty$ uniformly distributed if for any interval $J \subset I$
$$\lim_{N \to \infty} \frac{\#\{x_n \in J : n \leq N\}}{N} = |J|,$$
where $|J|$ is the Lebesgue measure of an interval $J$.

Uniformly distributed sequece $\{x_n\}_{n=1}^\infty$ is called of low discrepancy if
$$D_N = \sup_J \left| \frac{\#\{x_n \in J : n \leq N\}}{N} - |J| \right| = O\left(\frac{\log N}{N}\right).$$

It is well known that this is best possible and low discrepancy sequences play important role in numerical integrations.

We can construct a low discrepancy sequence using dynamical system.

**Definition 2.1.** The set $\{wx\}_{w \in \mathcal{W}}$ for which $wx$ exists and arranged in the order of words is called a van der Corput sequence defined by $F$.

When $F(x) = 2x \pmod 1$, $\mathcal{A}$ has two elements so we express $\mathcal{A} = \{0, 1\}$, and every words are admissible. Thus our van der Corput sequence is
$$x, 0x, 1x, 00x, 10x, 01x, 11x, 000x, \ldots,$$
and original van der Corput sequence is the case when $x = \frac{1}{2}$.

Our first goal is the following theorem:

**Theorem 2.2 ([4]).** *Let $F$ be a Markov and topologically transitive transformation with the same slope $|F'| \equiv \beta > 1$. Then, for any $x \in I$, the van der Corput sequence is of low discrepancy if and only if there exists no eigenvalue of the Perron–Frobenius operator in the annuls $e^{-\xi} < |z| \leq 1$ except $z = 1$.*

For Markov $\beta$–transformations, Ninomiya ([9]) proved the necessary and sufficient condition that the van der Corput sequence for $x = 0$ is of low discrepancy by direct calculation. He extended the results to non Markov $\beta$–transformations ([10]).

Rough sketch of the proof is the following. For any interval $J$
$$\begin{aligned} P^n 1_J(x) &= \sum_{y : F^n(y) = x} 1_J(y) |F^{n'}(y)|^{-1} \\ &= \beta^{-n} \sum_{|w|=n} 1_J(wx). \end{aligned}$$

Note that $\sum_{|w|=n} 1_J(wx)$ is the number of hits to $J$ in the van der Corput sequence corresponding to words with length $n$. On the other hand, as a rough expression,

$$P^n 1_J(x) = |J|\rho(x) + \eta^n,$$

where $\eta$ is the second greatest eigenvalue in modulus. Thus the number of $wx$ which hits to $J$ with $|w| = n$ equals $\beta^n |J|\rho(x) + (\beta\eta)^n$. The discrepancy depends on the term $(\beta\eta)^n$, and we already knew that $|\eta| \geq \beta^{-1}$. Thus $\eta = \beta^{-1}$ is the best possible case, and at that time, the number of $wx$ which hit to $J$ with $|w| \leq n$ equals

$$\frac{\beta^{n+1} - 1}{\beta - 1}|J|\rho(x) + n.$$

The number $N$ of admissible $wx$ with $|w| \leq n$ is of order $\beta^n$. This says the discrepancy is of order $\frac{\log N}{N}$.

## 3. Fredholm matrix

To make the proof of Theorem 2.2 rigorous, we define generating functions. Let $F$ be a piecewise linear Markov transformation with the same slope $\beta > 1$. Let for an interval $J$

$$s^J(z,x) = \sum_{n=0}^{\infty} z^n P^n 1_J(x).$$

Let $s(z,x)$ be a vector with coefficients $s^{\langle a \rangle}(z,x)$ ($a \in \mathcal{A}$). Then we get a renewal equation of the form:

$$s(z,x) = (I - \Phi(z))^{-1}\chi(x),$$

where $\Phi(z)$ is a $\mathcal{A} \times \mathcal{A}$ matrix and $\chi(x)$ is a $\mathcal{A}$ dimensional vector whose coefficients equal:

$$\Phi(z)_{a,b} = \begin{cases} z\beta^{-1} & \text{if } F(\langle a \rangle) \supset \langle b \rangle^o, \\ 0 & \text{otherwise}, \end{cases}$$

$$\chi(x)_a = 1_{\langle a \rangle}(x).$$

We call the matrix $\Phi(z)$ the Fredholm matrix associated with $F$.

**Theorem 3.1.**
$$\det(I - \Phi(z)) = \frac{1}{\zeta(z)}.$$

This theorem says that the eigenvalues of $P$ is determined by the zeros of $\det(I - \Phi(z))$. We can express $s^J(z,x)$ using the Fredholm matrix and a row vector whose coefficients are the functions with the radius of convergence $\beta$ ([4]).

From these, we can get the rigorous proof of Theorem 2.2.

We can extend the above discussion to more general cases for 1-dimensional transformations.

We use the signed symbolic dynamics defined by the orbits of the endpoints of $\langle a \rangle$ ($a \in \mathcal{A}$).

Let for $x \in I$

$$x^+ = \lim_{y \uparrow x} a_1^y a_2^y \cdots,$$
$$x^- = \lim_{y \downarrow x} a_1^y a_2^y \cdots.$$

For an interval $J$, we denote $(\sup J)^+$ and $(\inf J)^-$ by $J^+$ and $J^-$, respectively. Especially, for $a \in \mathcal{A}$, we denote by $a^+$ and $a^-$ instead of $\langle a \rangle^+$ and $\langle a \rangle^-$.

Let $\tilde{\mathcal{A}}$ be the set of $a^\sigma$ ($a \in \mathcal{A}, \sigma = \pm$).

We also define

$$s^{y^\sigma}(z,x) = \sum_{n=0}^{\infty} z^n \beta^{-n} \sum_{w \in \mathcal{W}_n} \sigma(y^\sigma, wx) \delta[\langle w[1] \rangle \supset \langle a_1^y \rangle, \exists \theta wx],$$

where $\theta$ is the shift and

$$\delta[L] = \begin{cases} 1 & \text{if } L \text{ is true,} \\ 0 & \text{if } L \text{ is false,} \end{cases}$$

$$\sigma(y^\sigma, x) = \begin{cases} +\frac{1}{2} & \text{if } y \geq_\sigma x, \\ -\frac{1}{2} & \text{if } y <_\sigma x, \end{cases}$$

$$x <_\sigma y = \begin{cases} x < y & \sigma = +, \\ x > y & \sigma = -. \end{cases}$$

Then

$$s^J(z,x) = s^{J^+}(z,x) + s^{J^-}(z,x).$$

We can also get the similar equation

$$s(z,x) = (I - \Phi(z))^{-1} \chi(z,x),$$

where

$$s(z,x) = (s^{a^\sigma}(z,x))_{a \in \mathcal{A}, \sigma = \pm},$$
$$\chi(z,x) = (\chi^{a^\sigma}(z,x))_{a \in \mathcal{A}, \sigma = \pm},$$
$$\chi^{y^\sigma}(z,x) = \begin{cases} \sigma(y^\sigma, x) & \text{if } \theta y^\sigma = (a_2^y)^{\sigma \operatorname{sgn} a_1^y} \\ \sum_{n=0}^{\infty} z^n \beta^{-n} \sigma(\theta^n y^\sigma, x) & \text{otherwise,} \end{cases}$$
$$\Phi(z) = (\phi_{a^\sigma, b^\tau}(z))_{a,b \in \mathcal{A}, \sigma, \tau = \pm},$$
$$\phi_{y^\sigma, b^\tau}(z) = \begin{cases} +z\beta^{-1}/2 & \text{if } b^\tau \leq_\sigma \theta y^\sigma, \text{ and } \theta y^\sigma = (a_2^y)^{\sigma \operatorname{sgn} a_1^y} \\ -z\beta^{-1}/2 & \text{if } b^\tau >_\sigma \theta y^\sigma, \text{ and } \theta y^\sigma = (a_2^y)^{\sigma \operatorname{sgn} a_1^y} \\ \sum_{n=1}^{\infty} z^n \beta^{-n} \sigma(\theta^n y^\sigma, b^\tau)(\operatorname{sgn} y^\sigma[1, n-1]) & \text{otherwise.} \end{cases}$$

Moreover, Theorem 3.1 also holds (cf. [3]). This leads to the general theory concerning the discrepancy of sequences generated by one-dimensional piecewise linear transformations.

**Definition 3.2.** We call an endpoint of $\langle a \rangle$ a Markov endpoint if the image of this point by $F^n$ for some $n \geq 1$ coincides with some endpoint of $\langle b \rangle$ ($b \in \mathcal{A}$).

**Theorem 3.3 ([2]).** *Let us denote by $k$ the number of non-Markov endpoints of $F$.*

1. *Let $\zeta_n \beta^{-n}$ be the n-th coefficient of $(1-z)\zeta(z)$. Assume that $\zeta_n$ is bounded, then*

$$D_N = O\left(\frac{(\log N)^{k+2}}{N}\right).$$

2. Let $\Phi_{11}(z)$ be the minor of the Fredholm matrix $\Phi(z)$ corresponding to Markov endpoints. Assume that $\det(I - \Phi_{11}(z)) \neq 0$ in $|z| < \beta$ and there exists no sigularity of $\zeta(z)$ in $|z| < \beta$ except 1. Then the discrepancy satisfies

$$D_N = O\left(\frac{(\log N)^{k+1}}{N}\right).$$

## 4. Higher-dimensional cases

We want to extend the result to higher-dimensional cases.

Let $I = [0,1]^d$. We will call a product of intervals also by intervals. We can define the Perron–Frobenius operator $P$ for higher-dimensional cases. Here we use Jacobian $\det F$ instead of derivative $F'$. For Markov cases, we can find the spectra of the Perron–Frobenius operator by the Fredholm matrix by almost the same discussion. Even for non Markov piecewise linear cases, instead of signed symbolic dynamics using the idea of screens we can construct Fredholm matrix and we can determine the spectra of the Perron–Frobenius operator ([6]).

We consider a space of functions for which there exists a decomposition $f = \sum_w C_w 1_{\langle w \rangle}$ such that for any $0 < r < 1$

$$\sum_w |C_w| r^{|w|} < \infty.$$

This is a slight extention of the space of functions with bounded variation in 1-dimensional cases. Here, note that, though the corresponding zeta function has meromorphic extention to $|z| < e^\xi$, even for the simplest Bernoulli transformation $F(x,y) = (2x, 2y) \pmod 1$, the essential spectrum radius is greater than $e^{-\xi}$.

Here is a question.

> Are there any transformation whose spectrum radius equals $e^{-\xi}$?

If such a transformation for which the Jacobian is constant exists, we will be able to construct van der Corput sequences of low discrepancy. Here we call a sequence of low discrepancy if

$$\sup_J \left| \frac{\#\{x_i \in J : i \leq N\}}{N} - |J| \right| = O\left(\frac{(\log N)^d}{N}\right),$$

where sup is taken over all intervals $J \subset I$.

It is proved that this is best possible for $d = 2$, and this will be also true for $d \geq 3$.

We have not yet the answer, but we have constructed low discrepancy sequences for $d = 2$ using dynamical system. Also for $d = 3$ we have constructed low discrepancy sequences using a sequence of transformations $F_1, F_2, \ldots$ but we have not yet constructed a transformation for which $F_n = F^n$.

We have not suceeded to get the answer to our question, but these results suggests that the question will be positive.

### 4.1. Two-dimensional cases

To construct low discrepancy sequences, we need a transformation not only expanding but also shuffling. For two-dimensional cases, let $s_0$ be the infinite sequence of 0's, and

$$s_1 = w_1 w_2 \cdots,$$

where $w_n$ is the word with length $2^{n-1}$ and only the last symbol equals 1 and other symbols are 0, that is,

$$s_1 = 101000100000001\cdots.$$

We consider the digitwise sum modulo 2 on the set of infinite sequences of 0 and 1. Let $\theta$ be a shift map to left. Then

$$s_0, s_1, \theta s_1, \ldots, \theta^{n-1} s_1$$

generate all the words with length $n$. Let us define

$$F\begin{pmatrix}x\\y\end{pmatrix} = \begin{pmatrix}\theta x\\\theta y\end{pmatrix} + \begin{pmatrix}s_{y_1}\\s_{x_1}\end{pmatrix},$$

where we identify $x \in [0,1)$ and its binary expansion $x_1 x_2 \cdots$. Then we can prove for $I = [\alpha, \alpha + 2^{-n+m}] \times [\beta, \beta + 2^{-n-m})$ for binary rationals $\alpha$ and $\beta$ and $m \le n$, the image of $I$ by $F^n$ does not overlap and $F^n(I) = [0,1)^2$. For example, $I = [0,1) \times [0, 1/4)$ ($n = m = 1$), though $\theta y$ belongs only to $[0, 1/2)$, $x_1$ takes both 0 and 1. Hence, $\theta y$ expands all $[0,1)$ by adding $s_0$ or $s_1$ depending on where $x$ belongs. Thus $F(I) = [0,1)^2$.

**Theorem 4.1 ([7]).** *The van der Corput sequence generated by the above $F$ is of low discrepancy.*

### 4.2. Three-dimensional cases

For three-dimensional cases, we have not yet proved to get low discrepancy sequences by one transformation. So we use a sequece of transformations $F, F_2, F_3, \ldots$. Even for higher dimensional cases, we will be able to construct low discrepancy sequences in the same way. However, we have no proof yet.

We denote

$$F_n \begin{pmatrix}x\\y\\z\end{pmatrix} = \begin{pmatrix}x'\\y'\\z'\end{pmatrix}.$$

Then expressing $x', y', z'$ in binary expansions, we define

$$\begin{pmatrix}x'_1\\y'_1\\z'_1\\x'_2\\y'_2\\z'_2\\\vdots\end{pmatrix} = \begin{pmatrix}x_{n+1}\\y_{n+1}\\z_{n+1}\\x_{n+2}\\y_{n+2}\\z_{n+2}\\\vdots\end{pmatrix} + M \begin{pmatrix}x_n\\y_n\\z_n\\\vdots\\x_1\\y_1\\z_1\\0\\\vdots\end{pmatrix}.$$

Here, $M$ is an infinite-dimensional matrix. One example of $M$ is following.

|  | $x_{-1}$ | $y_{-1}$ | $z_{-1}$ | $x_{-2}$ | $y_{-2}$ | $z_{-2}$ | $x_{-3}$ | $y_{-3}$ | $z_{-3}$ | $x_{-4}$ | $y_{-4}$ | $z_{-4}$ |
|---|---|---|---|---|---|---|---|---|---|---|---|---|
| $x_1$ | 0 | 1 | 1 | 0 | 1 | 0 | 0 | 1 | 0 | 0 | 0 | 1 |
| $y_1$ | 1 | 0 | 1 | 0 | 0 | 1 | 0 | 0 | 1 | 1 | 0 | 0 |
| $z_1$ | 1 | 1 | 0 | 1 | 0 | 0 | 1 | 0 | 0 | 0 | 1 | 0 |
| $x_2$ | 0 | 1 | 0 | 0 | 0 | 1 | 0 | 1 | 0 | 0 | 1 | 0 |
| $y_2$ | 0 | 0 | 1 | 1 | 0 | 0 | 0 | 0 | 1 | 0 | 0 | 1 |
| $z_2$ | 1 | 0 | 0 | 0 | 1 | 0 | 1 | 0 | 0 | 1 | 0 | 0 |
| $x_3$ | 0 | 1 | 0 | 0 | 1 | 0 | 0 | 0 | 1 | 0 | 1 | 0 |
| $y_3$ | 0 | 0 | 1 | 0 | 0 | 1 | 1 | 0 | 0 | 0 | 0 | 1 |
| $z_3$ | 1 | 0 | 0 | 1 | 0 | 0 | 0 | 1 | 0 | 1 | 0 | 0 |
| $x_4$ | 0 | 0 | 0 | 0 | 0 | 1 | 0 | 1 | 0 | 0 | 0 | 1 |
| $y_4$ | 0 | 0 | 0 | 1 | 0 | 0 | 0 | 0 | 1 | 1 | 0 | 0 |
| $z_4$ | 0 | 0 | 0 | 0 | 1 | 0 | 1 | 0 | 0 | 0 | 1 | 0 |
| $x_5$ | 0 | 0 | 0 | 0 | 1 | 0 | 0 | 1 | 0 | 0 | 1 | 1 |
| $y_5$ | 0 | 0 | 0 | 0 | 0 | 1 | 0 | 0 | 1 | 1 | 0 | 1 |
| $z_5$ | 0 | 0 | 0 | 1 | 0 | 0 | 1 | 0 | 0 | 1 | 1 | 0 |
| $x_6$ | 0 | 0 | 0 | 0 | 0 | 0 | 0 | 1 | 0 | 0 | 0 | 1 |
| $y_6$ | 0 | 0 | 0 | 0 | 0 | 0 | 0 | 0 | 1 | 1 | 0 | 0 |
| $z_6$ | 0 | 0 | 0 | 0 | 0 | 0 | 1 | 0 | 0 | 0 | 1 | 0 |
| $x_7$ | 0 | 0 | 0 | 0 | 0 | 0 | 0 | 1 | 0 | 0 | 1 | 0 |
| $y_7$ | 0 | 0 | 0 | 0 | 0 | 0 | 0 | 0 | 1 | 0 | 0 | 1 |
| $z_7$ | 0 | 0 | 0 | 0 | 0 | 0 | 1 | 0 | 0 | 1 | 0 | 0 |
| $x_8$ | 0 | 0 | 0 | 0 | 0 | 0 | 0 | 0 | 0 | 0 | 1 | 0 |
| $y_8$ | 0 | 0 | 0 | 0 | 0 | 0 | 0 | 0 | 0 | 0 | 0 | 1 |
| $z_8$ | 0 | 0 | 0 | 0 | 0 | 0 | 0 | 0 | 0 | 1 | 0 | 0 |
| $x_9$ | 0 | 0 | 0 | 0 | 0 | 0 | 0 | 0 | 0 | 0 | 0 | 1 |
| $y_9$ | 0 | 0 | 0 | 0 | 0 | 0 | 0 | 0 | 0 | 1 | 0 | 0 |
| $z_9$ | 0 | 0 | 0 | 0 | 0 | 0 | 0 | 0 | 0 | 0 | 1 | 0 |
| $x_{10}$ | 0 | 0 | 0 | 0 | 0 | 0 | 0 | 0 | 0 | 0 | 0 | 0 |
| $y_{10}$ | 0 | 0 | 0 | 0 | 0 | 0 | 0 | 0 | 0 | 0 | 0 | 0 |
| $z_{10}$ | 0 | 0 | 0 | 0 | 0 | 0 | 0 | 0 | 0 | 0 | 0 | 0 |

This $M$ has following properties: for any nonnegative integers $n$ and $m$: the determinant of the minor matrix with coordinates

$$\{x_1, \ldots, x_{n+m}\} \times \{y_{-1}, \ldots, y_{-n}, z_{-1} \ldots, z_{-m}\}$$

and the determinant of the minor matrix with coordinates

$$\{x_1, \ldots, x_n, y_1, \ldots, y_m\} \times \{z_{-1}, \ldots, z_{-n-m}\}$$

do not vanish. $M$ also has symmetry in the permutations $x$, $y$ and $z$. This matrix corresponds to $s_0$ and $s_1$ in two-dimensional case, and this shuffles coordinates, and we get for $k \geq n + m$

$$I = [\alpha, \alpha + 2^{-k-n-m}) \times [\beta, \beta + 2^{-k+n}) \times [\gamma, \gamma + 2^{-k+m})$$

or

$$I = [\alpha, \alpha + 2^{-k-n}) \times [\beta, \beta + 2^{-k-m}) \times [\gamma, \gamma + 2^{-k+n+m}),$$

then

$$F_k(I) = [0, 1)^3,$$

where $\alpha$, $\beta$ and $\gamma$ are binary rationals.

**Theorem 4.2.** *The van der Corput sequence generated by $F_1, F_2, \ldots$ is of low discrepancy.*